\documentclass[11pt,leqno]{article}  
  
\usepackage{amssymb}  
\usepackage{latexsym} 
\usepackage[notcite,notref]{showkeys} 
  
\newtheorem{theorem}{Theorem}[section]  
\newtheorem{lemma}[theorem]{Lemma}  
\newtheorem{corollary}[theorem]{Corollary}    
\newenvironment{proof}{\noindent {\bf  
Proof}.\ }{\proofbox}  
\newenvironment{proofof}[1]{\noindent {\bf  
Proof of #1}.}{\proofbox\par\smallskip\par}  
  
\def\numberlikeadb{\global\def\theequation{\thesection.\arabic{equation}}}  
\numberlikeadb

\newcommand{\halmos}{\rule{1ex}{1.4ex}}  
\newcommand{\proofbox}{\hspace*{\fill}\mbox{$\halmos$}}

\newcommand{\E}{\mathbb{E}}  
\newcommand{\pr}{\mathbb{P}}  
\def\re{\mathbb{R}} 
\newcommand{\m}[1]{\marginpar{\tiny{#1}}}

\newcommand{\ints}{\mathbb{Z}}  

\newcommand{\dtv}{d_{\mbox{{\rm\tiny TV}}}}

\newcommand{\eqa}{\begin{eqnarray}}  
\newcommand{\ena}{\end{eqnarray}}  
\newcommand{\eq}{\begin{equation}}  
\newcommand{\en}{\end{equation}}  
\newcommand{\eqs}{\begin{eqnarray*}}  
\newcommand{\ens}{\end{eqnarray*}}

\def\l{\lambda}  
\def\L{\Lambda}

\def\d{\delta}

\def\h{\eta}  
  
\def\th{\theta}  
  
\def\k{\kappa}  
\def\m{\mu}  
\def\n{\nu}  
\def\p{\pi}  
\def\r{\rho}  
\def\s{\sigma}  
\def\t{\tau}

\def\nin{\noindent}  
  
\def\msk{\medskip}  
  
\def\Blb{\left\{}  
\def\Brb{\right\}}  
  
\def\giv{\,|\,}  
\def\Giv{\,\Big|\,}  
\def\non{\nonumber}  
  
\def\Eq{\ =\ }  
\def\Le{\ \le\ }

\def\sji{\sum_{j\ge1}}

\def\Ref#1{{\rm (\ref{#1})}}  
  
\def\bp{\begin{proof}}  
\def\ep{\end{proof}}

\def\bone{{\bf 1}}  
\def\un{^{(n)}}

\def\ignore#1{}  

\def\Xm{X^\m} 
\def\Ge{\ \ge\ } 
\def\tso{\t_{\{s,0\}}} 
\def\tsm{\t_{\{s\}}^\m} 
\def\tss{\t_{\{s\}}}
 
\def\ftmj{\ff_{\tso^{\m,j-1}}} 
\def\xtm{{\Xm_t}} 
\def\xtms{{\Xm_{t - \tsm}}} 
\def\xtmv{{\Xm_{t-v}}} 
\def\ftsmj{g_{js}^\m} 
\def\law{{\cal L}} 
\def\ep{\hfil\halmos\break\vskip20pt} 
\def\ff{{\cal F}} 
\def\nno{N_0} 
\def\jj{{\cal J}} 
\def\nnc{N^n_c}
\def\nncc{(\nnc)\comp} 
\def\sjj{\sum_{j\in\jj}} 
\def\nnoc{\nno\cup\nnc}
\def\Z{{\mathbb Z}}
\def\P{{\mathbb P}}
\def\half{{\textstyle \frac12}}
\def\aaa{{\mathcal A}}
\def\smi{{\sum_{m\ge1}}}
\def\comp{^{{\rm c}}}
\def\hpo{{\widehat{\rm Po}}}
\def\remark{\nin {\bf Remark.\ }}
\def\nat{{\mathbb N}}
\def\qsd{quasi-stationary distribution}
\def\Def{\ :=\ }
\def\sjc{\sum_{j\in C}}
\def\skc{\sum_{k\in C}}
\def\lti{\lim_{t\to\infty}}
\def\tti{t\to\infty}
\def\etal{{\it et al.\/}}
\def\bir{b}
\def\dea{d}
\def\tK{{\widetilde K}}
\def\tk{{\tilde k}}
\def\klf{k_0}
\def\ds{_{\{s\}}}
\def\ui{^{(1)}}
\def\ut{^{(2)}}
\def\ex{\E}
\def\ul{^{(l)}}
\def\Bl{\left(}
\def\Br{\right)}
\def\alp{\alpha} 
\def\Aa{{\alp_+}}

\begin{document}  
  
\title{Total variation approximation for quasi-equilibrium distributions}  
\author{  
A. D. Barbour\footnote{Angewandte Mathematik, Universit\"at Z\"urich,  
Winterthurertrasse 190, CH-8057 Z\"URICH;  
ADB was supported in part by Schweizerischer Nationalfonds Projekt Nr.\  
20--107935/1.\msk}  
\ and  
P. K. Pollett\footnote{University of Queensland;  
PKP was supported in part by the Australian Research Council Centre of Excellence
for Mathematics and Statistics of Complex Systems.
}\\  
Universit\"at Z\"urich and University of Queensland  
}  
  
\date{}  
\maketitle  
  
\begin{abstract}  
Quasi--stationary distributions, as discussed by Darroch \& Seneta (1965),
have been used in biology to describe the steady state behaviour of population
models which, while eventually certain to become extinct, nevertheless
maintain an apparent stochastic equilibrium for long periods.  These
distributions have some drawbacks: they need not exist, nor be unique,
and their calculation can present problems.
In this paper, we give biologically plausible conditions under which
the \qsd\ is unique, and can be closely approximated by distributions
that are simple to compute.
\end{abstract}  
  
{\it Running head:}  Quasi-stationary distributions

\section{Introduction}  
\setcounter{equation}{0}  
The logistic growth model of Verhulst~(1838) was the first to describe
mathematically the evolution of a population to a non-zero equilibrium,
contrasting with the Malthusian law of exponential growth.
Its stochastic version, a Markov chain~$X$ in continuous time in which~$X(t)$
represents the number of individuals at time~$t$ in a population in a 
prescribed area~$A$, has transition rates
\eq
\begin{array}{cll}
    &q_{i,i+1} \Eq bi,\qquad q_{i,i-1} \Eq di + ei^2 /A,\qquad\qquad &i\ge1;\\
    &q_{ij} \Eq 0 &\mbox{otherwise},
\end{array}
\label{transitions}
\en
where $b$ and~$d$ are the {\em per capita\/} rates of birth and natural
mortality, and there is an additional {\em per capita\/} death rate~$ex$,
due to crowding, at population density~$x=i/A$. The stochastic model
has the drawback that its equilibrium distribution assigns probability~$1$ to
the state zero, population extinction, irrespective of the initial state. 
This apparently negates the most valuable property of Verhulst's model,
its ability to allow an equilibrium other than extinction.  
However, if $b>d$ and~$A$ is large, the population density~$X(t)/A$ can be 
expected to remain near the `carrying capacity' $\k := (b-d)/e$ for a very long time, 
in an apparent (and often biologically relevant) non-extinct stochastic equilibrium.

Darroch and Seneta~(1965), building on the work of Yaglom~(1947) in the context
of branching processes, introduced the concept of a
quasi-stationary distribution, in an attempt to reconcile
these at first sight inconsistent properties of the model.  In a discrete
time Markov chain~$X$
consisting of an absorbing state~$0$ together with a single finite transient
aperiodic class~$C$, the limiting conditional probabilities
\eq\label{q-limits}
   q_j \Def \lti \pr_i[X(t) = j \giv X(t) \in C], \quad i,j\in C,
\en
exist, and are the same for each~$i\in C$.  The~$q_j$, $j\in C$, also determine
a \qsd, in the sense that
\eq\label{qse}
   q_k \Eq \sjc q_j p_{jk} \Big/ \sjc q_j \skc p_{jk},
\en
where~$P := (p_{jk})$ denotes the one step probability transition matrix.
If, however, $C$ is countably infinite, the situation is very much less
satisfactory; there may be no \qsd, or exactly
one, or infinitely many, and determining which of these is the case
may be a difficult problem.  Even when there is a unique \qsd, its
calculation can pose substantial problems, unless the equations~\Ref{qse}
happen to have an obvious solution, because the probabilistic 
definition~\Ref{q-limits} involves conditioning on an event which, in
the limit as $\tti$, has probability zero.  This appears to make the 
\qsd\ unsatisfactory for typical biological 
applications.

In this paper, we give conditions, simply expressed in terms of the 
properties of the process~$X$, under which things are in fact
much simpler.  Under the conditions of Theorem~\ref{quasi-approx}, there is exactly
one \qsd, and it can be approximated to a specified accuracy by the equilibrium
distribution~$\p^\m$ of a `returned process'~$X^\m$.  What is more, under slightly
more stringent conditions, the distribution of~$X(t)$ is shown in 
Theorem~\ref{quasi-longtime} to be close to the \qsd\ 
for long periods of time. 

The returned process, introduced by Bartlett~(1960, pp.24-25) and used
by Ewens~(1963, 1964) in a population genetical setting, is a Markov process
that evolves exactly like~$X$,
up to the time at which~$0$ is reached, but is then instantly returned to
a random state in~$C$, chosen according to the probability measure~$\m$.
The mapping $\m\mapsto\p^\m$, studied in the paper of Ferrari \etal~(1995),
is contractive under our conditions, and iterating the mapping leads to the
unique \qsd\ $m$ on~$C$, which satisfies $m=\p^m$.  In many practical
applications, including the stochastic logistic model of~\Ref{transitions}
when~$A$ is large, iteration is unnecessary, inasmuch as {\em any\/}
distribution~$\p^\m$ is extremely close to~$m$. Furthermore, since~$\p^\m$ is
a genuine equilibrium distribution, its computation does not involve 
conditioning on sets of vanishing probability, and is hence typically
much simpler.

The main results, Theorems \ref{quasi-approx} and~\ref{quasi-longtime}, are proved in 
Section~\ref{return-measure}.    In Section~\ref{b&d}, as an illustration, we discuss 
the application of the theorems to birth and death processes, of which the stochastic
logistic model~\Ref{transitions} is an example.  These processes have the advantage
of having been widely studied, because of their relatively simple structure, allowing
our results to be easily interpreted; however, the theorems are equally applicable 
to processes with more complicated structure.

\section{A general approximation}\label{return-measure}  
\setcounter{equation}{0}  
Let~$X$ be a stable, conservative and non-explosive pure jump Markov process  
on a countable state space, consisting of a single transient class~$C$ together  
with a cemetery state~$0$. For any probability distribution~$\m$ on~$C$, define 
the modified process~$\Xm$ with state space~$C$ to have exactly the  
same behaviour as~$X$ while in~$C$, but, on reaching~$0$, to be instantly returned  
to~$C$ according to the distribution~$\m$.  Thus, if~$Q$ denotes the infinitesimal matrix 
associated with~$X$, and~$Q^\m$ that belonging to~$X^\m$, we have 
\eq\label{Qmu-def} 
   q_{ij}^\m \Eq q_{ij} + q_{i0}\m_j \qquad \mbox{for } i,j \in C. 
\en 
In this section, under a rather simple set of conditions,  we show that the 
quasi-stationary distribution~$m$ of~$X$ is unique, and
can be approximated in total variation to a prescribed accuracy by the {\em stationary\/}  
distribution of~$X^\m$, for an arbitrary choice of~$\m$.  We give a bound on the 
total variation distance between $m$ and~$\p^\m$ that is expressed solely in terms 
of hitting probabilities and mean hitting times for the process~$X$, and which is the 
same for all~$\m$.  The bound is such that it can be expected to be small in circumstances 
in which the process~$X$ typically spends a long time in~$C$ in apparent equilibrium, 
before being absorbed in~$0$ as a result of an `exceptional' event. If the bound is
not, as it stands, small enough for practical use, it can be improved geometrically 
fast by iteration of the return mapping $\m\mapsto\p^\m$.

Our basic conditions are as follows.
 
\medskip 
\nin{\bf Condition A}.\ \,{\sl There exist $s\in C$, $p>0$ and $T<\infty$ such that, 
uniformly for all $k\in C$,} 
\eqs 
     &{\rm (i)}&  p_k \Def \pr_k[X \mbox{ hits } s \mbox{ before } 0] \Ge p \,;\\ 
     &{\rm (ii)}&  \E_k[\tso]  \Le T \ <\ \infty . 
\ens

\nin Here,  $\pr_k$ and $\E_k$ refer to the distribution of~$X$ conditional on $X(0)=k$, and 
\eq\label{tau-def}
  \t_A \ :=\ \inf\{t>0\colon\,X(t) \in A, \,X(s) \notin A \mbox{ for some } s < t\}, 
\en 
the infimum over the empty set being taken to be~$\infty$.  Condition A~(i) can be expected
to be satisfied in reasonable generality; Condition A~(ii), although satisfied by the stochastic  
logistic model, is not so immediately natural.   

We now introduce the quantity
\eq\label{U-def} 
   U\ :=\ \sum_{k\in C} q_{k0} / \{q_k \E_k(\t_{\{k,0\}})\}. 
\en 
To interpret the meaning of~$U$, observe that a renewal argument for~$\Xm$, 
with renewal epochs the visits to any specific~$j\in C$, shows that 
\eq\label{equilibrium-bnd}
   \p^\m(j) \Le 1/\{q_j\E_j(\t_{\{j,0\}})\}. 
\en
In particular, if~$X$ has a \qsd~$m$, it follows from~\Ref{U-def} that
\[
   U \ \ge\ \sum_{i\in C} \p^m(i) q_{i0} \Eq \sum_{i\in C} m(i) q_{i0} \Eq \l_m,
\]
where~$\l_m$ is the rate at which the $X$-process, starting in the \qsd~$m$, leaves~$C$:
$\pr_m[X(t) \in C] = e^{-\l_m t}$.  Thus~$U$ acts as a computable upper bound for any~$\l_m$.
Note that $p,T$ and~$U$ are all quantities that can reasonably be bounded
using a knowledge of the process~$X$.

In the remainder of this section, we show that the \qsd~$m$ exists, 
is close to any~$\p^\m$, and well
describes the long time behaviour of~$X$ prior to absorption in~$0$, as long
as $UT/p$ is small enough. Our first  main result is the following.

\begin{theorem}\label{quasi-approx} 
Suppose that Condition~A is satisfied, and that
 $2UT/p < 1$.  Then $X$ has a unique \qsd~$m$, and,
for any probability measure~$\m$ on~$C$, we have 
\[ 
   \dtv(m,\p^\m) \Le 2UT/p.
\] 
\end{theorem} 

\remark Of course, for the theorem to imply that~$\p^\m$ is a sharp approximation
to~$m$, one needs~$U$  to be small enough (and therefore certainly finite).   
In many applications, $X$ can only jump to~$0$ 
from a small number of states in~$C$, and, if the quasi-equilibrium really 
behaves like a genuine equilibrium for long periods of time, 
the quantity $\E_k(\t_{\{k,0\}})$, for each such~$k$, can be expected to  
contain a large contribution 
from paths that, after leaving~$k$, spend a very long time `in equilibrium' in  
other states of~$C$, before either returning to~$k$ or being absorbed in~$0$.  
In such applications, as in the next section, these two features combine 
to make~$U$ small. 

To prove the theorem, we first need some preparatory results.
We first show that, under Condition~A, the mean time to hitting the state~$s$ 
is uniformly bounded,
for all return processes~$X^\m$, and for all initial states.

\begin{lemma}\label{s-mean} 
Under Condition~A, 
for all probability measures~$\m$ on~$C$ and for all $r\in C$, we have 
\[ 
    \E_r \t_{\{s\}}^\m \Le T/p \ <\ \infty, 
\]  
where $\t_A^\m$ is defined similarly to~$\t_A$, but with the process~$\Xm$ in place of~$X$. 
\end{lemma} 
 
\proof 
Recursively define 
\eqs 
   \tso^{\m,1} &:=& \tso^\m;\\ 
   \tso^{\m,j} &:=& \inf\Bigl\{ t>\tso^{\m,j-1}\colon\,\Xm(t) \in \{s,0\}, \,\Xm(u) \notin \{s,0\}  
     \mbox{ for some } \tso^{\m,j-1} < u < t \Bigr\},\\ 
   &&\hspace{2in}\mbox{}       \qquad j\ge2, 
\ens 
and, for $j\ge1$, let $Z_j := I[\Xm(\tso^{\m,l}) = 0, 1 \le l \le j]$, taking $Z_0=1$.   
Then it follows that 
\[ 
   \tss^\m \Eq \sji (\tso^{\m,j}  - \tso^{\m,j-1}) Z_{j-1}. 
\] 
Now $\E_r \tso^{\m,1} \le T$, by Condition~A(ii), and, for $j\ge2$, 
\eqs 
   \E\{(\tso^{\m,j}  - \tso^{\m,j-1}) Z_{j-1} \giv \ftmj\} 
       &=& Z_{j-1} \sum_{k\in C}\m_k \E_k\tso \Le TZ_{j-1}, 
\ens
by Condition~A(ii), where $\ftmj$ denotes the $\s$-field of events up to the 
stopping time $\tso^{\m,j-1}$.  Then, for $j\ge1$, 
\[ 
  \E\{Z_j \giv \ftmj\} \Le (1-p) Z_{j-1}, 
\] 
by Condition~A(i).  Hence, for $j\ge1$, and for any~$r\in C$, it follows that 
\[ 
   \E_r\{(\tso^{\m,j}  - \tso^{\m,j-1}) Z_{j-1}\} \Le T(1-p)^{j-1}, 
\] 
and so $\E_r \tss^{\m} \le T/p$, as required. 
\ep 
 
It follows in particular from Lemma~\ref{s-mean} that 
$\E_s \tss^\m \le T/p < \infty$, so that~$\Xm$ is 
positive recurrent on~$C$; denote its stationary distribution 
by~$\p^\m$.  Then, for any $f\colon\,C\to\re$ bounded, set 
\[ 
  h_f^\m(j) \ :=\ -\int_0^\infty \{\E_j f(X_t^\m) - \p^\m(f)\}\,dt, 
  \qquad j \in C.
\] 
To show that the integral is well-defined, note that 
\[ 
   |h_f^\m(j)| \Le \int_0^\infty 2\|f\|\, \dtv(\law(\Xm(t)\giv \Xm(0)=j),\p^\m)\,dt,
\] 
where $\|\cdot\|$ denotes the supremum norm.
The latter integral is finite provided that $\E_s\{(\tss^\m)^2\} < \infty$,
by the coupling inequality (Lindvall~2002, (2.8)) and from Pitman~(1974, Corollary~1,
(1.23) with $r=2$).  That this is the case follows from the next lemma.

\begin{lemma}\label{s-second} 
Under Condition~A, 
for all probability measures~$\m$ on~$C$, we have 
\[ 
    \E_s\{ (\tss^\m)^2\} \ <\ \infty. 
\]  
\end{lemma} 

\proof
Writing $\t := \tss^\m$, note that 
\eq\label{tau-squared}
   \mbox{}\quad\t^2 \Eq \Bl \int_0^\infty I[\t > t]\,dt \Br^2
    \Eq 2\int_0^\infty I[\t > t] \Bl \int_t^\infty I[\t > u]\,du \Br\,dt.
\en
Now, from Lemma~\ref{s-mean} and by the Markov property, we have 

\[
   \E \Bl \int_t^\infty I[\t > u]\,du \Giv \ff^\m_t \Br \Le (T/p) I[\t > t],
\]
where $\ff^\m_t$ denotes the history of~$X^\m$ up to time~$t$.
Hence, taking expectations in~\Ref{tau-squared}, it follows that
\[
   \E_s(\t^2) \Le 2(T/p)\E_s\t \Le 2(T/p)^2,
\]
again from Lemma~\ref{s-mean}, completing the proof.
\ep

\nin It is shown in the proof of Theorem~\ref{quasi-longtime} that
the distribution of~$\tss^\m$ actually has an exponential tail.

The functions~$h_f^\m$ are central to the argument to come.
First, we show that they are bounded and Lipschitz, with
appropriate constants. 
 
\begin{lemma}\label{h-diff} 
For all $j\in C$, 
\[ 
  |h_f^\m(j) - h_f^\m(s)| \Le 2\|f\|\,T/p. 
\] 
\end{lemma} 
 
\proof 
For any $j\in C$, we can write 
\eqa 
  -h_f^\m(j) &=& \int_0^\infty \E_j \{(f(\xtm)- \p^\m(f))\,I[\tsm \le t]\}\,dt \non\\ 
     &&\qquad\mbox{} + \int_0^\infty \E_j \{(f(\xtm)- \p^\m(f))\,I[\tsm > t]\}\,dt. 
    \label{h-split} 
\ena 
Then, by the strong Markov property, we have 
\eqs 
  \E_j \{(f(\xtm)- \p^\m(f))\,I[\tsm \le t]\} 
     &=& \E_j\Blb \E_s \{(f(\xtms)- \p^\m(f))\,I[\tsm \le t]\}\Brb \\ 
  &=& \int_0^t \ftsmj(v) \E_s(f(\xtmv) - \p^\m(f))\,dv, 
\ens 
where $\ftsmj$ denotes the probability density of the random variable~$\tsm$ 
for the process started at~$j$.  Hence it follows that 
\eqs 
  \lefteqn{\int_0^\infty \E_j \{(f(\xtm)- \p^\m(f))\,I[\tsm \le t]\}\,dt}\\ 
  &=& \int_0^\infty dt \int_0^\infty dv\, \ftsmj(v)\E_s(f(\xtmv) - \p^\m(f)) \bone\{v\le t\}. 
\ens 
Now, since 
\eqs 
  \lefteqn{\int_0^\infty \bone\{v\le t\}|\E_s(f(\xtmv) - \p^\m(f))|\,dt}\\ 
   &&\Le 2\|f\|\,\int_0^\infty \dtv(\law(\xtm\giv \Xm(0)=s),\p^\m) \ <\ \infty, 
\ens 
we can use Fubini's theorem to conclude that 
\eqs 
  \lefteqn{\int_0^\infty \E_j \{(f(\xtm)- \p^\m(f))\,I[\tsm \le t]\}\,dt}\\ 
  &=& \int_0^\infty  \ftsmj(v) \Blb\int_v^\infty  \E_s(f(\xtmv) - \p^\m(f))\,dt\Brb\,dv \\ 
  &=& \int_0^\infty  \ftsmj(v) h_f^\m(s)\,dv \Eq h_f^\m(s). 
\ens 
Hence, from~\Ref{h-split} and Lemma~\ref{s-mean}, it follows that 
\[ 
  |h_f^\m(j) - h_f^\m(s)| \Le 2\|f\|\, \E_j(\tsm) \Le 2\|f\|\,T/p, 
\] 
as required. 
\ep 
\nin In particular, the function~$h_f^\m$ is itself bounded. 
 
\medskip 
A similar argument, by conditioning on the time of the first jump, shows that 
\eq\label{h-equation} 
  h_f^\m(j) \Eq -q_j^{-1}\{f(j) - \p^\m(f)\} + \sum_{k\in C,k\ne j}q_j^{-1}q_{jk}^\m h_f^\m(k), 
\en 
where $q_j := \sum_{k\in C\cup\{0\}} q_{jk} < \infty$ because~$X$ is conservative, 
and the sum in~\Ref{h-equation} is absolutely convergent because~$h_f^\m$ is bounded. 
This can be rewritten in the form 
\eq\label{Stein-eq} 
  (Q^\m h_f^\m)(j) \Eq f(j) - \p^\m(f),\qquad j\in C, 
\en
so that, for any bounded~$f$ and for any probability measures $\m$ and~$\n$ on~$C$, we 
have  
\eq\label{mu-nu} 
   \p^\m(Q^\n h_f^\n) \Eq \p^\m(f) - \p^\n(f).
\en 
In the terminology of Stein's method, \Ref{Stein-eq} determines~$h_f^\m$ to
be the solution~$h$ of the Stein equation $(Q^\m h)(j) \Eq f(j) - \p^\m(f)$ for
the distribution~$\p^\m$, corresponding to the given function~$f$.
Also, by Dynkin's formula, we have 
\eq\label{Dynkin} 
   \p^\m(Q^\m h) \Eq 0 
\en 
for any bounded function~$h$ (for the special case $h=h_f^\m$,  
this follows from~\Ref{mu-nu}).  
These considerations put us into a position to prove Theorem~\ref{quasi-approx}. 
 
\medskip
\proofof{Theorem~\ref{quasi-approx}} 
Take any probability measures $\m$ and~$\n$ on~$C$.  Then~\Ref{mu-nu} gives 
\[ 
   \p^\m(Q^\n h_f^\n) \Eq \p^\m(f) - \p^\n(f), 
\] 
whereas~\Ref{Dynkin} gives $\p^\m(Q^\m h_f^\n) = 0$.  Taking the difference, 
we obtain 
\eq\label{pi-diff} 
  \p^\m(f) - \p^\n(f) \Eq \p^\m(Q^\n h_f^\n - Q^\m h_f^\n).  
\en 
Now, for bounded~$h$ and for any $i\in C$, 
\eqs 
   (Q^\n h - Q^\m h)(i) &=&  
    \sum_{k\in C} q_{ik}^\n (h(k) - h(i)) - \sum_{k\in C} q_{ik}^\m (h(k) - h(i)), 
\ens 
with both sums absolutely convergent, and, from~\Ref{Qmu-def}, it then follows that 
\eqs
   (Q^\n h - Q^\m h)(i) &=&  q_{i0}\sum_{k\in C}(\n(k)-\m(k))(h(k) - h(i))\\
    &=&  q_{i0}\sum_{k\in C}(\n(k)-\m(k))(h(k) - h(s)), 
\ens 
since $\skc \n(k) = \skc \m(k) = 1$.
Hence, from~\Ref{pi-diff}, we have 
\[ 
  \p^\m(f) - \p^\n(f) \Eq \sum_{i\in C} \p^\m(i) q_{i0}
      \sum_{k\in C}(\n(k)-\m(k))(h_f^\n(k) - h_f^\n(s)), 
\] 
and, from Lemma~\ref{h-diff}, this gives  
\eq\label{penultimate} 
  |\p^\m(f) - \p^\n(f)| \Le \sum_{i\in C} \p^\m(i) q_{i0}\, 2\|f\|(T/p)\,\|\n-\m\|_{TV}. 
\en 
Thus it follows that
\eq\label{contraction}
   \|\p^\n-\p^\m\|_{TV} \Le (2T/p)\sum_{i\in C} \p^\m(i) q_{i0}\,\|\n-\m\|_{TV},
\en
and~\Ref{equilibrium-bnd} then implies that 
\eq\label{final} 
    \|\p^\n-\p^\m\|_{TV} \Le (2UT/p)\,\|\n-\m\|_{TV}.
\en
This, by the Banach fixed point theorem, establishes the first part of the theorem,
and the second follows by taking $\n = m$, and using the fact that,
for probability measures $F$ and~$G$, $\dtv(F,G) = \half\|F-G\|_{TV}$.
\ep

\medskip
We now turn our attention to the distribution of~$X(t)$ for fixed values of~$t$,
starting from any initial distribution, and compare it to~$m$, the distribution 
at any time of the return process~$X^m$ started in the \qsd~$m$.
We begin by taking the initial state of~$X$ to be~$s$, and remark later that this 
restriction makes little difference, provided that~$s$ is hit at least once.  

\begin{theorem}\label{quasi-longtime}
Let $B := Tq_s/p \ge 1$.  Then,
 under Condition~A and if also $2UT/p < 1$, there is a constant~$K$ such that
\[
  \dtv(\law_s(X(t)),m) \Le 
        Ut + KB\sqrt{ \frac T{pt}} + (2/e)^{pt/16T} \ =:\ \h(t).
\]
\end{theorem}

\remark
Hence, if $UB^2T/p \ll 1$, the distribution $\law_s(X(t))$ is close to~$m$ 
for all times~$t$ such that
\[
    B^2 T/p\ \ll\ t\ \ll\ U^{-1}.
\]

\medskip
 
\proof 
The argument is based on coupling two copies $X\ui$ and~$X\ut$ of the return 
process~$X^m$, with~$X\ui$ in equilibrium and with~$X\ut$ starting in~$s$.  
The coupling is achieved by forcing~$X\ui$ to follow the same sequence of states 
as~$X\ut$ after the first time that it hits~$s$, and to have identical residence
times in all states other than~$s$. Define
\[
   \tss^0(1) \Def \inf\{t\ge0\colon\, X\ui(t) = s\};\qquad \tss^0(2) := 0, \Eq v,
\]
and let
\[
   \tss^n(l) \Def \inf\{t\ge\tss^{n-1}(l)\colon\, X\ul(t-) \neq X\ul(t) = s\},
            \quad l=1,2,
\]
denote the $n$-th return time of~$X\ul$ to~$s$.
Then, conditional on the event that $\tss^0(1) = v$, we have
\eqs
   \lefteqn{\dtv(\law(\tss^n(1)\giv \tss^0(1)=v),\law(\tss^n(2)))}\\
    &&\Eq \dtv(\d_v*q_s^{-1}G(n,1), q_s^{-1}G(n,1)) 
    \Le c_G q_s v n^{-1/2},
\ens
for a suitable constant~$c_G$,
where $G(n,1)$ is the Gamma distribution with shape parameter~$n$ and unit scale
parameter, and $\d_v$ is the point mass at~$v$.  Hence, for any $n\ge1$,
we can couple $X\ui$ and~$X\ut$ by arranging that $\tss^n(1) = \tss^n(2)$,
with the two processes to be run identically thereafter, and the probability of
this coupling failing, conditional on $\tss^0(1) = v$, is at most
$c_G q_s v n^{-1/2}$.  Thus, in particular,
\eq\label{1st-bnd}
  \phantom{HH} \dtv(\law(X\ui(t)),\law(X\ut(t))) \Le \pr[\tss^n(2) > t] 
       + c_Gq_s n^{-1/2} T/p,
\en
using Lemma~\ref{s-mean}.  It now remains to show that we can reach the bound given
in the theorem by choosing $n$ almost as a multiple of~$t$.

Now $\tss^n(2)$ is a sum of independent random variables, each with distribution
$\law_s(\t^m\ds)$, where $\t^m\ds$ is defined as in Lemma~\ref{s-mean}.
By that lemma and Markov's inequality, it follows that
\[
   \pr_r[\t^m\ds \ge 2T/p] \ \le\ 1/2,\qquad r\in C,
\]
and hence that
\[
    \pr_s[\t^m\ds \ge 2kT/p] \Le 2^{-k} 
       \Eq \exp\Blb -\frac{p\log2}{2T}\Bl \frac{2kT}p \Br\Brb \qquad\mbox{for all}\ k\ge1.
\]
Thus the distribution~$\law_s(\t^m\ds)$ is stochastically bounded above by that
of 
\[
     \frac{2T}p \Blb 1 + \frac1{\log 2}E \Brb,
\]
where~$E$ has a standard exponential distribution.  Hence the distribution of~$\tss^n(2)$
is stochastically bounded above by that of
\[
    \frac{2T}p \Blb n + \frac1{\log 2}G_n \Brb, 
\]
where $G_n \sim G(n,1)$.  The inequality $\pr[G_n \ge 2n] \le (2/e)^n$ thus implies
that
\eq\label{gamma-bnd}
   \pr\left[\tss^n(2) > \frac{2Tn}p \Blb 1 + \frac2{\log 2} \Brb \right]
     \le (2/e)^n.
\en

So, for any fixed~$t$, using $1 + 2/\log 2 \le 4$, we take $n = n_t 
:= \lfloor tp/8T \rfloor$ in~\Ref{1st-bnd}, giving
\eq\label{2nd-bnd}
   \dtv(\law(X\ui(t)),\law(X\ut(t))) \Le (2/e)^{n_t} 
       + c_Gq_s n_t^{-1/2} T/p,
\en
from which it follows that, for $t \ge 16T/p$,
\eq\label{3rd-bnd}
   \dtv(\law(X\ui(t)),\law(X\ut(t))) \Le (2/e)^{pt/16T} 
       + 4c_G\,\frac{(Tq_s/p)^{3/2}}{\sqrt{q_st}}.
\en

We first observe that $\law(X\ui(t)) = m$ for all~$t$.
Then we have
\[
   \pr[\t_{\{0\}}(1) \le t] \Eq 1 - e^{-\l_m t} \Le Ut,
\]
where $\t_{\{0\}}(1) := \inf\{t\ge0\colon\,X\ui(t)=0\}$.
On the event that $X\ui$ and~$X\ut$ are successfully coupled at
$\tss^{n_t} \le t$, it thus follows that the event that 
neither hits~$0$ before~$t$ has probability at least $1-Ut$, 
and, on this event, $X\ut(t)$ is also the value of an $X$-process
starting in~$s$, since~$X\ut$ has had no visits to~$0$ before~$t$.
This,
\ignore{
whereas 
\[
   \dtv(\law(X\ut(t)),\law_s(X(t))) = \pr_s[\t_{\{0\}} \le t].
\]
However, each excursion starting from~$s$ begins with a residence
time in~$s$, and these are independent of each other and of the 
sequence of states visited, being exponentially distributed with 
mean~$q_s^{-1}$. Hence, setting $r := \lceil 2q_st \rceil$,  we have
\[
    \pr_s[\tss^r(2) \le t] \Le \pr[G_r \le q_st] \Le e^{-q_st/4},
\]
and
\[
    \pr_s[\t_{\{0\}} \le \tss^r(2)] \Le 2UTr,
\]
by Lemma~\ref{return-prob} and because $2UT/p < 1$.  
These observations,
} 
together with~\Ref{3rd-bnd}, completes the proof.
\ep

\remark
Denoting by $A(\{s\},\{0\})$ the event that~$X$ hits $s$ before~$0$,
the same argument can be used to show that
$\dtv(\law_k(X(t) \giv A(\{s\},\{0\})),\law_s(X^m(t)))$ is at most~$\h(t)$
for any $k\in C$, under the conditions of Theorem~\ref{quasi-longtime}.  
Hence, conditional on the event that $X$ hits~$s$
before reaching~$0$,  the distribution of~$X(t)$ starting from any $k\in C$
is also close to~$m$ for all times~$t$ such that
\[
    B^2 T/p\ \ll\ t\ \ll\ U^{-1},
\]
provided that $UB^2T/p \ll 1$.  Thus the \qsd~$m$ is then indeed the
appropriate long time approximation to the distribution of~$X$ in~$C$,
for times $t \ll U^{-1}$.

Note also that the coupling used in Theorem~\ref{quasi-longtime} may be
very pessimistic, only making use of the residence times in~$s$.
For most processes, the variability in the remaining residence times
and in the possible sequences of states can be exploited to get sharper
bounds.  However, in the examples for which we make computations below,
the quantity $B^2 T/p$ is of only polynomial order in the size of the system,
whereas $U^{-1}$ is exponentially large; hence even this crude estimate is
more than adequate.

\section{Birth and death processes}\label{b&d} 
\setcounter{equation}{0} 

Consider now a birth and death process with $C=\{1,2,\dots,N\}$ (for
$C=\nat$, replace $N$ by $\infty$ in what follows) having
birth rates $\bir_j > 0$, $1\le j < N$, with $\bir_0 = 0$ and
$\bir_N=0$ if $N<\infty$, and with strictly positive death rates $\dea_j$,
$j\in C$. It is convenient to introduce the quantities 
$(\alp_j,\ j\in C)$, where $\alp_1=1$ and, for $j>1$,
$$
\alp_j = \frac{\bir_1 \cdots \bir_{j-1}}{\dea_2 \cdots \dea_j}.
$$
The return process with $\m = \d_{\{1\}}$, equivalent to 
re-defining~$d_1$ to be zero, is then recurrent if $\Aa := \sji \alp_j < \infty$,
in which case $\p^\m(j) \Eq \alp_j/\Aa$, so that its computation is very
easy. We now wish to investigate when this distribution can be used as
a reasonable approximation to the effective steady state behaviour of the process. 

In order to apply Theorems \ref{quasi-approx} and~\ref{quasi-longtime},
we need to choose a state~$s\ge1$, and find values for $p,T,B$ and~$U$.
For~$p$, let $r_k$, $k\ge1$, be the probability that the process starting in~$k$ 
hits $s$ before it hits $0$, where $s\geq1$. If $k> s$, then
$r_k=1$. Otherwise, $r_0=0$, $r_s=1$ and 
$$
  (\bir_k + \dea_k) r_k \Eq \bir_k r_{k+1} + \dea_k r_{k-1},\qquad k=1,2,\dots,s-1, 
$$
leading to $r_k = \sigma_k/\sigma_s$, where
$$
  \sigma_0 \Eq 0 \quad\mbox{and}\quad \sigma_k 
  \Eq \sum_{j=1}^{k} \frac{1}{\dea_j \alp_j}
     \quad \mbox{for}\ k=1,\dots,s. 
$$
Since $\sigma_k$ is non-decreasing in~$k$, we can take
\eq\label{p-bd}
  p \Eq r_1 \Eq 1/(\dea_1\sigma_s), 
\en
for any state $s\in C$.

For~$T$, 
we note that
\eq\label{T1-def}
  T_1 \Def \max_{1\le k\le s} \E_k\{\t_{\{s,0\}}\}
   \Le \sum_{j=1}^{s-1}\frac1{\bir_i\alp_i}\sum_{j=0}^i \alp_j \ <\ \infty,
\en
and that, for $k>s$,
$$
\E_{k} (\tau_{\{s,0\}})
  \Eq \E_{k} (\tau_{\{s\}})
  \Eq \sum_{j=s+1}^{k} \E_{j} (\tau_{\{j-1\}})
  \Eq \sum_{j=s+1}^{k} \frac{1}{\dea_{j} \alp_{j}} \sum_{i=j}^\infty \alp_i
$$
(Anderson~1991, Chapter~8). Since the latter quantity is
increasing in $k$, we may take
\begin{equation}
  T \Def \max(T_1,T_2), \quad\mbox{where}\quad
     T_2 \Def \sum_{j=s+1}^{\infty} \frac{1}{\dea_{j} \alp_{j}} \sum_{i=j}^\infty \alp_i.
     \label{pkp2}
\end{equation}
Note that then Condition~A(ii) holds if~$T_2$, the so-called ``$D$ series'',
converges, and that $T_2<\infty$ is a necessary and sufficient
condition for a birth and death process to have a unique 
\qsd\ (van Doorn~1991, part~2 of Theorem~3.2).  Note also that~$B := Tq_s/p$
can be bounded using \Ref{p-bd}--\Ref{pkp2}, together with the fact that
$q_s = \bir_s + \dea_s$.

Finally, the quantity $U$ can be evaluated as
\eq\label{U-bd}
U \Eq \frac{\dea_1}{ (\bir_1 + \dea_1) \,\E_{1} (\tau_{\{1,0\}}) }
  \Eq \frac{\dea_1}{ 1+ \bir_1 \,\E_{2} (\tau_{\{1\}}) }
  \Eq \frac{\dea_1}{ \sum_{j=1}^N \alp_j},
\en
because, also from Anderson~(1991, Chapter~8),
$$
\E_{i} (\tau_{\{i-1\}}) 
  \Eq \frac{1}{\dea_i \alp_i} \sum_{j=i}^N \alp_j
  \Eq \frac{1}{\bir_{i-1} \alp_{i-1}} \sum_{j=i}^N \alp_j,
$$
and in particular, since $\alp_1=1$,
$$
   1+ \bir_1 \E_{2} (\tau_{\{1\}}) \Eq 1 + \sum_{j=2}^N \alp_j 
  \Eq \sum_{j=1}^N \alp_j.
$$

In order to apply Theorems \ref{quasi-approx} and~\ref{quasi-longtime}
in practice, we need to be able to
bound the quantities $p,T,B$ and~$U$ by assigning concrete expressions 
in terms of the $\bir_j$ and~$\dea_j$ to replace \Ref{p-bd}--\Ref{U-bd}.  
Simple estimates can be derived under the
assumptions that the death rates $\dea_j$ are increasing in~$j$, and
that the ratios $\bir_j/\dea_j$ are decreasing, with $\bir_1/\dea_1 > 1$. 
If this is the case,
define $s\ge1$ in such a way that \hbox{$\bir_s/\dea_s \ge 1 > \bir_{s+1}/\dea_{s+1}$,}
and let $1\le s_1 \le s < s_2$ be such that 
$$
  \bir_{s_1}/\dea_{s_1} \ =:\ \r_1 \ >\ 1 \ >\ \r_2 \Def \bir_{s_2}/\dea_{s_2}.
$$
Then $x_j := \dea_{j+1}\alp_{j+1}/\dea_1$ is maximal at $j=s$, and
$$
  x_j \ \ge\ \r_1^{s_1},\quad s_1 \le j \le s;\qquad x_j/x_s \Le \r_2^{j-s_2},
  \quad j\ge s_2.
$$
Hence, from \Ref{p-bd}--\Ref{U-bd}, we have the bounds
\eqs
  p &=& 1\Big/ \sum_{j=0}^{s-1} x_j^{-1} \Ge \frac{\r_1}{\r_1-1}
     \Blb 1 + (s-s_1)\r_1^{-s_1} \Brb^{-1};\\
  U &\le& \dea_{s_1}(\r_1-1)\r_1^{-s_1};\\
  T_2 &\le& \{s_2-s+1/(1-\r_2)\} \sum_{j\ge s+1} \frac1{\dea_j};\\
  T_1 &\le& \sum_{j=1}^{s-1}\frac1{\bir_i\alp_i}\sum_{j=0}^i \alp_j
    \Le \{s - s_1 + \r_1/(\r_1-1)\}\sum_{i=1}^{s-1}\frac1{\bir_i}.
\ens
Thus if, for instance, $\dea_j$ grows at most polynomially fast in~$j$,
with the sum $\sji \dea_j^{-1} < \infty$, and if $s_1$ and $s-s_1$ are large and
of comparable size, then $T/p$ is roughly of order~$s^2$ and $B=Tq_s/p$ of polynomial 
order in~$s$, whereas~$U$
is geometrically small with~$s$, making $UT/p$ very small indeed.

More precise calculations for the stochastic logistic model 
of~\Ref{transitions} as $A\to\infty$, with $s = \lfloor \k A \rfloor$, give
$$
  T \Eq O(\log A),\quad p \ge 1 - d/b,\quad 
      U \Le \Blb 1 + \frac{b-d}{b+d}\Brb^{-\k A/2},\ 
     B \Eq O(A\log A),
$$
so that $UT/p$ is geometrically small in~$A$ as $A\to\infty$. Thus,
for the stochastic logistic model, the unique \qsd\ can be very closely
approximated by any return distribution, as long as~$A$ is large.
Entirely similar estimates are true for the SIS epidemic model,
which models the number of susceptibles in a closed population
of size~$N$, to be thought of as large but finite.  The process is 
a birth and death process on $\{0,1,\ldots,N\}$ having rates
\eq\label{SIS}
   \bir_i \Def \l i(1 - i/N)\quad\mbox{and}\quad \dea_i \Def \m i, \qquad 0 \le i \le N;
\en
in this case, $UT/p$ is geometrically small in~$N$ if $\m < \l$, and there
is a \qsd\ close to $s := \lfloor N(1-\m/\l) \rfloor$.

\section*{Acknowledgement}  
ADB wishes to thank the School of Mathematical Sciences at Monash University,
the Australian Research Council Centre of Excellence
for Mathematics and Statistics of Complex Systems,
and the Institute Mittag--Leffler for their warm hospitality, while part of this 
work was accomplished.

\end{document}